\documentclass[12pt]{amsart}
\usepackage{amsmath,amsthm}
\usepackage{amssymb}
\newtheorem{theorem}{Theorem}
\newtheorem{corollary}[theorem]{Corollary}

\usepackage{amssymb}
\usepackage{amsmath}
\usepackage{breqn}
\usepackage[a4paper, nohead, headsep=0.6cm, left = 2cm, right = 2cm, top = 2cm, bottom = 2cm, twoside]{geometry}
\usepackage{inputenc}
\usepackage[english]{babel}

\begin{document}
\title{INTERPOLATION PROPERTIES OF CERTAIN CLASSES OF NET SPACES}

\author{A. H. KALIDOLDAY, E. D. NURSULTANOV}
\address{L.N.Gumilyov Eurasian National University, Nur-Sultan, Kazakhstan} 
\address{M.V.Lomonosov Moscow State University (Kazakhstan branch), Nur-Sultan, Kazakhstan}

\email{aitolkynnur@gmail.com, er-nurs@yandex.ru}
\thanks{This research was supported by the Ministry of Education and Science of the Republic of Kazakhstan (projects no. AP08956157 and AP08856479).}

\renewcommand{\thefootnote}{}



\renewcommand{\thefootnote}{\arabic{footnote}}
\setcounter{footnote}{0}

\maketitle

\begin{abstract}
The paper studies the interpolation properties of net spaces $N_{p,q}(M)$, when $M$ is the set of dyadic cubes in $\mathbb{R}^n$, and also when $M$ is the family of all cubes with parallel faces to the coordinate axes in $\mathbb{R}^n$.
It is shown that, in the case when $M$ is the set of dyadic cubes the scale of spaces is closed with respect to the real interpolation method. In the case, when $M$ is the set of all cubes with parallel faces to the coordinate axes, an analogue of the Marcinkiewicz-Calderon theorem on cones of non-negative functions is given.

\end{abstract}

\section{Introduction}

Let in $\mathbb{R}^n$  is given $n$-dimensional Lebesgue measure  $\mu,$ $M$ is the set of all cubes in $\mathbb{R}^n.$ Further $M$ will be called as a {\itshape "net".} For function $f(x),$ defined and integrable on each $e$ from $M,$ we define the function

$$
\Bar{f}(t,M)=\sup_{ \begin{subarray}{1} e\in M \\ |e|\ge t \end{subarray} }\frac{1}{|e|} \bigg| \int_{e} f(x) dx \bigg|, \;\;\;\;  t>0,
$$
where the supremum is taken over all $e\in M$,  whose measure is $|e| \stackrel{\mathrm{def}}{=} \mu e \ge t$. In the case when $\sup \{|e| : e \in M\}=\alpha< \infty$ and $t>\alpha$ assuming that $\Bar{f}(t,M)=0.$ 

Let $p,q$ parameters satisfy the conditions $0<p\le \infty$, $0<q\le \infty$. Let's define the net spaces $N_{p,q}(M)$, as a set of all functions $f$, such that for $q<\infty$

$$
\|f\|_{N_{p,q}(M)}=\bigg( \int_0^{\infty} \bigg( t^{\frac{1}{p}} \Bar{f}(t,M) \bigg)^q \frac{dt}{t} \bigg)^{\frac{1}{q}}<\infty,
$$
and for $q=\infty$
$$
\|f\|_{N_{p,\infty}(M)}=\sup_{t>0} t^{\frac{1}{p}} \Bar{f}(t,M) < \infty.
$$

These spaces were introduced in the work  \cite{N1}.

Net spaces have found important applications in various problems of harmonic analysis, operator theory and the theory of stochastic processes \cite{AR1, AR2, ARN, NK, NTl, N2, N3, NT, NA}.

In this paper, we study the interpolation properties of these spaces. It should be noted here, that net spaces are in a sense close to the Morrey space:
$$
M_p^\alpha = \left\{ f: \;\;\; \sup_{y\in \mathbb{R}^n, \;\; t>0} t^{-\lambda} \left( \int_{|x+y| \le t} |f(x)|^p dx  \right)^{\frac{1}{p}} < \infty \right\} .
$$

In the case when $f(x)\ge 0,$ for $\frac{1}{p}= 1- \frac{\lambda}{n}$ 
$$
\|f\|_{N_{p,\infty}(M)} \asymp \|f\|_{M_1^\lambda}.
$$

The question of interpolation of Morrey spaces was considered in the works \cite{St, CM, P, RV, BRV, LR2}. It follows from the results of  \cite{P} that
$$
(M_{p}^{\lambda_0}, M_{p}^{\lambda_1})_{\theta,\infty}\hookrightarrow M_{p}^{\lambda}\,,
$$
where $\lambda=(1-\theta)\lambda_0+\theta \lambda_1$. In the works \cite{RV,BRV}  it was established that this inclusion is strict.

For net spaces $N_{p,q}(M)$, where $M$  is an arbitrary system of measurable sets from $\mathbb{R}^2$, we also have an embedding (see \cite[Theorem 1]{N1})

\begin{equation}\label{eq1}
    \left( N_{p_0,q_0}(M), N_{p_1,q_1}(M) \right)_{\theta,q} \hookrightarrow N_{p,q}(M),
\end{equation}
where $\frac{1}{p}=\frac{1-\theta}{p_0} + \frac{\theta}{p_1},$ $0<\theta<1,$   $0<q\le \infty$.

From (\ref{eq1}) it follows that if the linear operator $T$ bounded from  $A_i$ to $N_{p_i,\infty}(M),$ $i=0,1,$ then the operator $T$ bounded from $A_{\theta,q}$ to $N_{p,q}(M),$ where $\frac{1}{p}=\frac{1-\theta}{p_0} + \frac{\theta}{p_1}.$

The question arises whether the following equality will take place
\begin{equation}\label{eq2}
    \left( N_{p_0,q_0}(M), N_{p_1,q_1}(M) \right)_{\theta,q} = N_{p,q}(M).
\end{equation}

Here, in contrast to the Morrey spaces in the one-dimensional case, when $M$ is the set of all segments, the answer is positive  \cite{NT}.  

In this paper we show that, if  $M$ is the set of dyadic cubes in $\mathbb{R}^n$,  then the relation (\ref{eq2}) holds. In case, if $M$ is the set of all cubes, an analogue of the Marcinkiewicz-Calderon theorem on the cones of non-negative functions is obtained.

Given functions $F$ and $G$, in this paper $F \lesssim G$ means that $F \le CG$, where $C$ is a positive
number, depending only on numerical parameters, that may be different on different occasions. Moreover, $F \asymp G$ means that $F \lesssim G$ and $G \lesssim F.$

\section{Main results}

Let $(A_0, A_1)$ -- compatible pair of Banach spaces \cite{BL}.
$$
K(t,a;A_0,A_1)=\inf_{a=a_0+a_1}(\|a_0\|_{A_0}+t\|a_1\|_{A_1}), \  a\in A_0+A_1,
$$
-- functional Petre.
For $0<q<\infty,$  $0<\theta<1$

$$
(A_0,A_1)_{\theta,q}=\bigg\{ a\in A_0+A_1 : \|a\|_{(A_0,A_1)_{\theta,q}}=\bigg( \int_0^{\infty} (t^{-\theta} K(t,a))^q \frac{dt}{t}\bigg)^{1/q} < \infty \bigg\},
$$
and for $q=\infty$
$$
(A_0,A_1)_{\theta,q}=\bigg\{ a\in A_0+A_1 : \|a\|_{(A_0,A_1)_{\theta,q}}= \sup_{0<t<\infty} t^{-\theta} K(t,a) < \infty \bigg\}.
$$

A family of sets from $\mathbb{R}^n$ of the form 
$$
Q^m_{k}= \left[ \frac{k_1}{2^m}, \frac{k_1+1}{2^m}  \right) \times ... \times \left[ \frac{k_n}{2^m}, \frac{k_n+1}{2^m} \right),
$$
where $k \in \mathbb{Z}^n$, $m\in \mathbb{Z}$, is called the family of dyadic cubes and denoted by $M.$

Note that for arbitrary $m\in \mathbb{Z}$ the space $\mathbb{R}^n$  can be represented in the form
$$
\mathbb{R}^n=\bigcup_{k\in \mathbb{Z}^n}Q^m_{k},
$$
and the measure of intersection
$$\left| Q^m_{k} \cap Q^m_{r} \right|= \begin{cases} 2^{nm}, \ k_i=r_i, \ i=\overline{1,n} \\ 0, \  \text{in other cases}.\end{cases} 
$$

Let $m\in \mathbb{Z}$, cubes $Q^m_{k}$, $k\in \mathbb{Z}^n$ are called cubes of $m$  order.

Note also that if $n\ge m$, then each cube of $n$ order  is partitioned into  $4^{n-m}$  cubes of $m$ order.

\begin{theorem}\label{th1} Let $0<p_0<p_1< \infty $ and $0<q_0, q_1, q \le \infty.$ Let $M$ be the family of dyadic cubes.
Then
$$
\big( N_{p_0, q_0}(M), N_{p_1, q_1} (M)\big)_{\theta, q} = N_{p,q}(M),
$$
where $\frac{1}{p}= \frac{1-\theta}{p_0}+\frac{\theta}{p_1}, \ \theta \in(0,1)$.

\end{theorem}

The following statement is an attempt to answer the question about the interpolation of net spaces, when $M$ -- is the family of all cubes with parallel faces to the coordinate axes in $\mathbb{R}^n$.

\begin{theorem}\label{th2}
Let $n \le p_0 <p_1 < \infty,$ $0<q \le \infty,$ $M$ be the family of all cubes with parallel faces to the coordinate axes in $\mathbb{R}^n$. Let $G=\{ f : f(x) \ge 0 \},$  then for any $f \in G \bigcap N_{p,q}(M)$  it is true 

\begin{equation}\label{th2eq1}
\|f\|_{(N_{p_0,q_0}(M), N_{p_1, q_1}(M))_{\theta,q}} \asymp \|f\|_{N_{p,q}(M)},
\end{equation}
where the corresponding constants depend only on $p_i, q_i, \theta, q, i=0,1.$
\end{theorem}

The following corollary holds from Theorem \ref{th2}.

\begin{corollary}
Let $n \le p_0 <p_1 < \infty,$ $1 \le q_0, q_1 < \infty,$ $q_0 \neq q_1,$ $0< \tau, \sigma< \infty.$  $M$ and $G$ -- sets from Theorem \ref{th2}. If the following inequalities hold for a quasilinear operator 

\begin{equation}\label{th3eq1}
   \|Tf\|_{N_{q_0, \infty}(M)} \le F_0 \|f\|_{N_{p_0, \sigma}(M)}, \  f\in N_{p_0, \sigma}(M), 
\end{equation}
\begin{equation}\label{th3eq2}
    \|Tf\|_{N_{q_1, \infty}(M)} \le F_1 \|f\|_{N_{p_1, \sigma}(M)}, \  f\in N_{p_1, \sigma}(M),
\end{equation}
then for any $f \in G \cap N_{p,\tau}$ we have
\begin{equation}\label{th3eq3}
    \|Tf\|_{N_{q, \tau}(M)} \le c F_0^{1-\theta} F_1^{\theta} \|f\|_{N_{p, \tau}(M)},
\end{equation}
where $\frac{1}{p}= \frac{1-\theta}{p_0}+\frac{\theta}{p_1}, \ \theta \in(0,1)$ and the corresponding constant depends only on $p_i, q_i, \sigma, i=0,1.$
\end{corollary}

\section{Proof of Theorem 1}
\begin{proof}
Let $M$ --  be the set of dyadic cubes, $1<p_0<\infty$. Let us prove first
\begin{equation}\label{eq11}
    \left( N_{1,\infty}(M), N_{\infty, \infty}(M) \right)_{\theta,q}=N_{p,q}(M),
\end{equation}
where $\frac{1}{p}= 1-\theta$, $\theta \in (0,1).$

Let $m\in \mathbb{Z}$,  then the Euclidean space $\mathbb{R}^n$ is partitioned into disjoint cubes of order $m$ from $M$
$$
\mathbb{R}^n=\bigcup_{k\in \mathbb{Z}^n}Q^m_{k}.
$$
Let $f \in  N_{p,q}(M)$, define the function 
$$
\varphi_0(x)=\frac{1}{|Q^m_{k}|} \int_{Q^m_{k}} f(x)dx, \;\;\ x\in Q^m_{k}, \;\;\; k\in \mathbb{Z}^n.
$$
Then taking into account that the measure  $|Q^m_{k}|=2^{nm}$ we have
\begin{equation}\label{proof1,1}
    |\varphi_0(x)| \le \Bar{f}(2^{nm},M), \;\;\;\;  x\in \mathbb{R}^n,
\end{equation}
and
 $$
  \int_{Q^m_{k}} (f(x)-\varphi_0(x))dx=0.
 $$
 
 For the Petre functional, we have the following
$$
  K(t,f; N_{p_0, \infty}(M), N_{\infty,\infty}(M)) = \inf_{f=f_0+f_1} \big( \|f_0\|_{N_{p_0,\infty}(M)} + t\|f_1\|_{N_{\infty,\infty}(M)} \big) \le
 $$
 $$
 \le  \sup_{s >0} s^{\frac{1}{p_0}} \overline{(f-\varphi_0)} (s,M) + t\sup_{s >0}  \Bar{\varphi_0}(s,M)= \sup_{s >0} s^{\frac{1}{p_0}} \overline{(f-\varphi_0)} (s,M)+ t \bar{f}(2^{nm},M).
 $$
Consider the first term
$$
 \sup_{s >0} s^{\frac{1}{p_0}} \overline{(f-\varphi_0)}(s,M) \asymp \sup_{2^{nm} \ge s >0} s^{\frac{1}{p_0}} \overline{(f-\varphi_0)}(s,M) + \sup_{s \ge 2^{nm}} s^{\frac{1}{p_0}} \overline{(f-\varphi_0)} (s,M).
$$

Let $I$-- be an arbitrary cube from $M$ such that $|I|\ge 2^{nm}.$  Hence, $I$ is some cube of order $n$, where $n\ge m$. Taking into account that each dyadic cube of $n$ order  is partitioned into mutually disjoint cubes of $m$ order,  we obtain
$$
\left| \int_I (f-\varphi_0)(x)dx \right|
 =\bigg| \sum_{Q^m_{k}\subset I} \int_{Q^m_{k}} (f(x)-\varphi_0(x))dx  \bigg|=0.
$$
 
Hence,
$$
\sup_{s >0} s^{\frac{1}{p_0}} \overline{(f-\varphi_0)}(s,M) =   \sup_{2^{nm} \ge s>0} s^{\frac{1}{p_0}}  (\overline{f-\varphi_0)}(s,M)
$$$$
\le \sup_{2^{nm} \ge s>0} s^{\frac{1}{p_0}}  \bar{f}(s,M)+  \sup_{2^{nm} \ge s>0} s^{\frac{1}{p_0}}  \bar{\varphi_0}(s,M) \le
$$$$
\le \sup_{2^{nm} \ge s>0} \frac{1}{p_0} \int_0^s t^{\frac{1}{p_0}}\bar{f}(t,M) \frac{dt}{t} +\sup_{2^{nm} \ge s>0} s^{\frac{1}{p_0}}  \bar{f}(2^{nm},M) 
$$$$
\le \sup_{2^{nm} \ge s>0} \frac{1}{p_0} \int_0^s t^{\frac{1}{p_0}}\bar{f}(t,M) \frac{dt}{t} + 2^{\frac{nm}{p_0}}\bar{f}(2^{nm},M)=
$$$$
=\frac{1}{p_0}\int_0^{2^{nm}} t^{\frac{1}{p_0}}\bar{f}(t,M) \frac{dt}{t} +\frac{1}{p_0}\int_0^{2^{nm}} t^{\frac{1}{p_0}}\bar{f}(t,M) \frac{dt}{t} = \frac{2}{p_0} \int_0^{2^{nm}} \Bar{f}(y,M) y^{\frac{1}{p_0}-1} dy.
$$

In this way,
$$
K(a^m,f; N_{p_0, \infty}(M), N_{\infty,\infty}(M)) \le  c \int_0^{2^{nm}} y^{\frac{1}{p_0}} \Bar{f}(y,M) \frac{dy}{y}  + a^m \bar{f}(2^{nm},M),
 $$
where $a=2^{ \frac{n}{p_0} }>1$. Then we have

$$
\|f\|_{(N_{p_0, \infty}(M), N_{\infty,\infty}(M))_{\theta,q}} \asymp \left( \sum_{m\in \mathbb{Z}}\left( a^{-\theta m} K(a^m,f)\right)^q \right)^{\frac{1}{q}}
$$

$$
\le \left( \sum_{m\in \mathbb{Z}} \left( a^{-\theta m}  \left( c \int_0^{4^m} y^{\frac{1}{p_0}} \Bar{f}(y,M) \frac{dy}{y}  + a^m \bar{f}(2^{nm},M) \right)  \right)^q \right)^{\frac{1}{q}}.
$$

Applying the Minkowski inequality, we obtain the following

$$
\|f\|_{(N_{p_0, \infty}(M), N_{\infty,\infty}(M))_{\theta,q}}\le \left( \sum_{m\in \mathbb{Z}} \left( a^{-\theta m}  c \int_0^{4^m} y^{\frac{1}{p_0}} \Bar{f}(y,M) \frac{dy}{y}  \right)^q \right)^{\frac{1}{q}}
$$$$
+\left( \sum_{m\in \mathbb{Z}} \left( a^{m(1-\theta)}\bar{f}(4^m,M)  \right)^q \right)^{\frac{1}{q}}.
$$

Further, taking into account that $a=2^{ \frac{n}{p_0}}$  and applying Hardy's inequality for the first term, we obtain
$$
\|f\|_{(N_{p_0, \infty}(M), N_{\infty,\infty}(M))_{\theta,q}}\le  \left( \sum_{m\in \mathbb{Z}} \left( a^{-\theta m}  \sum^m_{k=-\infty} 2^{\frac{nk}{p_0}} \Bar{f}(2^{nk},M)  \right)^q \right)^{\frac{1}{q}}
$$$$
+ \left( \sum_{m\in \mathbb{Z}} \left( a^{m(1-\theta)}\bar{f}(2^{nm},M)  \right)^q \right)^{\frac{1}{q}}
$$$$
\le c \left( \sum_{m\in \mathbb{Z}} \left( 2^{\left(\left(-\theta\right) \frac{1}{p_0}+\frac{1}{p_0} \right) nm} \Bar{f}(2^{nm},M)  \right)^q \right)^{\frac{1}{q}}+ \left( \sum_{m\in \mathbb{Z}} \left( 2^{\frac{nm(1-\theta)}{p_0}}\bar{f}(2^{nm},M)  \right)^q \right)^{\frac{1}{q}}
$$$$
=c \left( \sum_{m\in \mathbb{Z}} \left( 2^{\frac{nm}{p}} \Bar{f}(2^{nm},M)  \right)^q \right)^{\frac{1}{q}} \asymp \|f\|_{N_{p,q}},
$$
where $\frac{1}{p}= \frac{1-\theta}{p_0}, \ \theta \in(0,1)$.

So we got the embedding
$$
N_{p,q} \hookrightarrow (N_{1, \infty}(M), N_{\infty,\infty}(M))_{\theta,q},
$$
where $\frac{1}{p}= 1-\theta, \ \theta \in(0,1)$. 
The reverse embedding follows from (\ref{eq1}). Hence the relation (\ref{eq11}) holds. To prove the general case, we use the reiteration theorem \cite[Theorem 3.5.3]{BL}.

Let $1<p_0<p_1< \infty$.
 From (\ref{eq11}) it follows that there are $\theta_0, \theta_1 \in (0,1)$  such that 
\begin{equation}
\begin{split}
    (N_{1, \infty}(M), N_{\infty, \infty}(M))_{\theta_0, q_0} = N_{p_0,q_0}(M) \\
     (N_{1, \infty}(M), N_{\infty, \infty}(M))_{\theta_1, q_1} = N_{p_1,q_1}(M),
\end{split}
\end{equation}
then by the reiteration theorem it follows that
$$
(N_{p_0, q_0}(M), N_{p_1,q_1}(M))_{\theta, q} = (N_{1, \infty}(M), N_{\infty,\infty}(M))_{\eta, q}= N_{p,q}(M).
$$

In the last equality, we took into account that  $\eta= (1-\theta)\theta_0 + \theta \theta_1$.

\end{proof}

\section{Proof of Theorem 2}
\begin{proof}

Let $\tau>0$, euclidean space $\mathbb{R}^n=\bigcup_{k=1}^{\infty} I_k $ partitions into non-intersecting half-open cubes $\{ I_k \}^{\infty}_{k=1}$  with faces parallel to the coordinate axes and such that $ | I_k|= \tau$. Let $f \in G \bigcap N_{p,q}(M)$, define the function
$$
\varphi_0(x)=\frac{1}{|I_k|} \int_{I_k} f(x)dx, \ x\in I_k, \;\;\; k\in \mathbb{N}.
$$

Then it is obvious that
\begin{equation}\label{proof2,1}
   \varphi_0(x) \le \Bar{f}(\tau), \;\;\;\;  x\in \mathbb{R}^n,
\end{equation}
and
 $$
 \int_{I_k} (f(x)-\varphi_0(x))dx=0.
 $$
 
Let $0<\sigma<\min{\{q_0,q_1,q}\}$, then
$$
  K(t,f; N_{p_0, \sigma}(M), N_{p_1, \sigma}(M)) = \inf_{f=f_0+f_1} \big( \|f_0\|_{N_{p_0,\sigma}(M)} + t\|f_1\|_{N_{p_1,\sigma}(M)} \big) 
 $$$$
 \le \left( \int_0^\infty \left( s^{\frac{1}{p_0}} \overline{(f-\varphi_0)} (s)\right)^\sigma \frac{ds}{s}\right)^{\frac{1}{\sigma}} +  t \left( \int_0^\infty \left( s^{\frac{1}{p_1}} \Bar{\varphi}_0(s)\right)^\sigma \frac{ds}{s}\right)^{\frac{1}{\sigma}}.
 $$

Estimate the first term, taking into account the inequality (\ref{proof2,1}), we have
 $$
\left( \int_0^\infty \left( s^{\frac{1}{p_0}} \overline{(f-\varphi_0)} (s)\right)^\sigma \frac{ds}{s}\right)^{\frac{1}{\sigma}} \asymp \left( \int_0^\tau \left( s^{\frac{1}{p_0}} \overline{(f-\varphi_0)} (s)\right)^\sigma \frac{ds}{s}\right)^{\frac{1}{\sigma}} 
$$$$
+ \left( \int_\tau^\infty \left( s^{\frac{1}{p_0}} \overline{(f-\varphi_0)} (s)\right)^\sigma \frac{ds}{s}\right)^{\frac{1}{\sigma}} = \left( \int_0^\tau \left( s^{\frac{1}{p_0}} \bar{f} (s)\right)^\sigma \frac{ds}{s}\right)^{\frac{1}{\sigma}} 
$$$$
+ \left( \int_0^\tau \left( s^{\frac{1}{p_0}} \bar{\varphi}_0 (s)\right)^\sigma \frac{ds}{s}\right)^{\frac{1}{\sigma}} + \left( \int_\tau^\infty \left( s^{\frac{1}{p_0}} \overline{(f-\varphi_0)} (s)\right)^\sigma \frac{ds}{s}\right)^{\frac{1}{\sigma}}
$$$$
\le \left( \int_0^\tau \left( s^{\frac{1}{p_0}} \bar{f} (s)\right)^\sigma \frac{ds}{s}\right)^{\frac{1}{\sigma}} + \Bar{f}(\tau)  \left( \int_0^\tau  s^{\frac{\sigma}{p_0}-1}  ds \right)^{\frac{1}{\sigma}} 
$$$$
+ \left( \int_\tau^\infty \left( s^{\frac{1}{p_0}} \overline{(f-\varphi_0)} (s)\right)^\sigma \frac{ds}{s}\right)^{\frac{1}{\sigma}}
= 2 \left( \int_0^\tau \left( s^{\frac{1}{p_0}} \bar{f} (s)\right)^\sigma \frac{ds}{s}\right)^{\frac{1}{\sigma}} 
$$$$ +   \left( \int_\tau^\infty \left( s^{\frac{1}{p_0}} \overline{(f-\varphi_0)} (s)\right)^\sigma \frac{ds}{s}\right)^{\frac{1}{\sigma}}.
$$

Let $I$ be an arbitrary cube from $M$,   $\partial I$ -- its boundary, then
$$
\bigg| \int_I (f-\varphi_0)(x)dx \bigg|=
 $$$$
 =\bigg| \sum_{I_k\subset I} \int_{I_k} (f(x)-\varphi_0(x))dx + \sum_{|I_k\cap \partial I| \ne 0} \int_{I_k\cap I} (f(x)- \varphi_0(x))dx \bigg|.
$$

Note that the first sum is equal to zero, and the second contains at most  $2^{n}\left( \left( \frac{|I|}{\tau} \right)^{\frac{1}{n}}+1 \right)^{n-1}$  terms. Moreover,  taking into account the non-negativity of the function  $f$, the second sum is estimated as follows
$$
\bigg| \sum_{|I_k\cap \partial I| \ne 0}  \int_{I_k\cap I} (f(x)- \varphi_0(x))dx \bigg| 
$$

$$
\le \sum_{|I_k\cap \partial I| \ne 0} \int_{I_k\cap I} f(x)dx  + \sum_{|I_k\cap \partial I| \ne 0} |I_k\cap  I| \frac{1}{|I_k|}  \int_{I_k} f(x)dx  
$$

$$
\le 2 \sum_{|I_k\cap \partial I| \ne 0} \int_{I_k} f(x)dx  \le 2\sum_{|I_k\cap \partial I| \ne 0} \tau \Bar{f}(\tau)
\le 4^{n} |I|^{\frac{n-1}{n}} \tau^{\frac{1}{n}} \Bar{f}(\tau).
$$

Hence, 
$$
\sup_{|I| \ge s} \frac{1}{|I|} \bigg| \int_I (f-\varphi_0)(x)dx\bigg| \le  \sup_{|I| \ge s} \frac{1}{|I|} 4^{n} |I|^{\frac{n-1}{n}} \tau^{\frac{1}{n}}  \Bar{f}(\tau)
$$$$
= 4^{n} \tau^{\frac{1}{n}}  \Bar{f}(\tau) \sup_{|I| \ge s} \frac{1}{|I|^{\frac{1}{n}}} \le  4^{n} \cdot \frac{\tau^{\frac{1}{n}} }{s^{\frac{1}{n}}}\Bar{f}(\tau).
$$

 Further, taking into account that $p_0 \ge n$ we have
 $$
 \left( \int_\tau^\infty \left( s^{\frac{1}{p_0}} \overline{(f-\varphi_0)} (s)\right)^\sigma \frac{ds}{s}\right)^{\frac{1}{\sigma}}= 
 \left( \int_\tau^\infty \left( s^{\frac{1}{p_0}} \sup_{|I| \ge s} \frac{1}{|I|} \bigg| \int_I (f-\varphi_0)(x)dx\bigg| \right)^\sigma \frac{ds}{s}\right)^{\frac{1}{\sigma}}
$$

$$
\le 4^{n} \left( \int_\tau^\infty \left( s^{\frac{1}{p_0}} \frac{\tau^{\frac{1}{n}} }{s^{\frac{1}{n}}}\Bar{f}(\tau) \right)^{\sigma} \frac{ds}{s}\right)^{\frac{1}{\sigma}} = 4^{n} \tau^{\frac{1}{n}} \Bar{f}(\tau) \left( \int_{\tau}^{\infty} \left( s^{\frac{1}{p_0}-\frac{1}{n}}\right)^\sigma \frac{ds}{s}\right)^{\frac{1}{\sigma}}
$$

$$
\asymp  \tau^{\frac{1}{n}} \Bar{f}(\tau) \tau^{\frac{1}{p_0}-\frac{1}{n}}=  \tau^{\frac{1}{p_0}} \Bar{f}(\tau)= c \left( \int_0^{\tau} \left( s^{\frac{1}{p_0}} \bar{f} (s)\right)^{\sigma} \frac{ds}{s}\right)^{\frac{1}{\sigma}}.
$$

Hence,
$$
\left( \int_0^\infty \left( s^{\frac{1}{p_0}} \overline{(f-\varphi_0)} (s)\right)^\sigma \frac{ds}{s}\right)^{\frac{1}{\sigma}} \le c \left( \int_0^\tau \left( s^{\frac{1}{p_0}} \bar{f} (s)\right)^\sigma \frac{ds}{s}\right)^{\frac{1}{\sigma}},
$$
where $c$ depends only on parameters $p_0,\sigma$.

To estimate the second term, we first show that,
\begin{equation}\label{proof2}
    \Bar{\varphi_0}(s)\le \begin{cases} \Bar{f}(\tau), \ \text{at  }  s\le \tau, \\ 4\Bar{f}(s), \ \text{at  }  s > \tau.
    \end{cases}
\end{equation}

 For $ s \le \tau$  from (\ref{proof2,1}) we have $\Bar{\varphi_0}(s)\le \Bar{f}(\tau).$ Let $s > \tau$, $Q\in  M,$ $|Q|=s $ and by $Q'$  we denote a cube, whose center coincides with the center of $Q$ and has an edge twice as large, then we have
$$
\Bar{\varphi_0}(s)= \sup_{\begin{subarray}{l} |Q| \ge s \\ Q\in M \end{subarray} } \frac{1}{|Q|} \sum_{I_k \cap Q \neq 0} \int_{I_k \cap Q} \frac{1}{I_k} \int_{I_k} f(x) dxdy 
$$$$
\le \sup_{\begin{subarray}{l} |Q| \ge s \\ Q\in M \end{subarray} } \frac{1}{|Q|} \sum_{I_k \cap Q \neq 0} \int_{I_k} f(x) dx  \frac{|I_k \cap Q|}{I_k} 
\le \sup_{\begin{subarray}{l} |Q| \ge s \\ Q\in M \end{subarray} } \frac{1}{|Q|} \sum_{I_k \cap Q \neq 0} \int_{I_k} f(x) dx  
$$$$
\le \sup_{\begin{subarray}{l} |Q| \ge s \\ Q\in M \end{subarray} } \frac{1}{|Q|}\int_{Q'} f(x) dx \le 4 \Bar{f}(s).
$$

Further, using the relation (\ref{proof2}) we obtain

$$
\left( \int_0^\infty \left( s^{\frac{1}{p_1}} \Bar{\varphi}_0(s)\right)^\sigma \frac{ds}{s}\right)^{\frac{1}{\sigma}} = \left( \int_0^\tau \left( s^{\frac{1}{p_1}} \Bar{\varphi}_0(s)\right)^\sigma \frac{ds}{s}\right)^{\frac{1}{\sigma}} 
$$$$
+ \left( \int_\tau^\infty \left( s^{\frac{1}{p_1}} \Bar{\varphi}_0(s)\right)^\sigma \frac{ds}{s}\right)^{\frac{1}{\sigma}} \le \Bar{f}(\tau) \left( \int_0^\tau \left( s^{\frac{1}{p_1}} \right)^\sigma \frac{ds}{s}\right)^{\frac{1}{\sigma}} 
$$$$
+ 4 \left( \int_\tau^\infty \left( s^{\frac{1}{p_1}} \Bar{f}(s) \right)^\sigma \frac{ds}{s}\right)^{\frac{1}{\sigma}}
= c \tau^{\frac{1}{p_1}} \Bar{f}(\tau) + 4 \left( \int_\tau^\infty \left( s^{\frac{1}{p_1}} \Bar{f}(s) \right)^\sigma \frac{ds}{s}\right)^{\frac{1}{\sigma}},
$$
where c depends only on parameters $p_1, \sigma.$

Hence,
$$
K(t,f; N_{p_0, \sigma}(M), N_{p_1,\sigma}(M)) \le c \left( \int_0^\tau \left( s^{\frac{1}{p_0}} \bar{f} (s)\right)^\sigma \frac{ds}{s}\right)^{\frac{1}{\sigma}} 
$$$$
+ t \left(  c \tau^{\frac{1}{p_1}} \Bar{f}(\tau) + 4 \left( \int_\tau^\infty \left( s^{\frac{1}{p_1}} \Bar{f}(s) \right)^\sigma \frac{ds}{s}\right)^{\frac{1}{\sigma}}   \right).
 $$

Taking into account the monotonicity of $\Bar{f}(s)$, for $\tau = t^{\big(\frac{1}{p_0}-\frac{1}{p_1}\big)^{-1} }$  we obtain 
$$
\|f\|_{(N_{p_0, \sigma}(M), N_{p_1,\sigma}(M))_{\theta,q}}\le \bigg( \int_0^{\infty} \bigg[ t^{-\theta} \bigg( \bigg( \int_0^{t^{\big(\frac{1}{p_0}-\frac{1}{p_1}\big)^{-1} }} \left( s^{\frac{1}{p_0}} \bar{f} (s)\right)^\sigma \frac{ds}{s}\bigg)^{\frac{1}{\sigma}} 
$$

$$
+t \bigg(  c t^{\big(\frac{1}{p_0}-\frac{1}{p_1}\big)^{-1}\frac{1}{p_1}} \Bar{f}(t^{\big(\frac{1}{p_0}-\frac{1}{p_1}\big)^{-1}}) + 4 \left( \int_{t^{\big(\frac{1}{p_0}-\frac{1}{p_1}\big)^{-1}}}^\infty \left( s^{\frac{1}{p_1}} \Bar{f}(s) \right)^\sigma \frac{ds}{s}\right)^{\frac{1}{\sigma}}   \bigg) \bigg) \bigg]^q \frac{dt}{t}\bigg)^{\frac{1}{q}}.
$$

Making the replacement  $\gamma= t^{\big(\frac{1}{p_0}-\frac{1}{p_1}\big)^{-1} }$  and applying the Minkowski inequality, we arrive at the following
$$
\|f\|_{(N_{p_0, \sigma}(M), N_{p_1,\sigma}(M))_{\theta,q}} \le \left( \int_0^{\infty} \left( \gamma^{-\theta \left(\frac{1}{p_0}-\frac{1}{p_1}\right)} \bigg( \int_0^{\gamma} \left( s^{\frac{1}{p_0}} \bar{f} (s)\right)^\sigma \frac{ds}{s}\bigg)^{\frac{1}{\sigma}} \right)^q \frac{d\gamma}{\gamma}  \right)^{\frac{1}{q}}
$$$$
+ \left( \int_0^{\infty} \left( \gamma^{\frac{1}{p}} \bar{f}(\gamma) \right)^q \frac{d\gamma}{\gamma}  \right)^{\frac{1}{q}} +
\left( \int_0^{\infty} \left( \gamma^{(1-\theta) \left(\frac{1}{p_0}-\frac{1}{p_1}\right)} \bigg( \int_\gamma^{\infty} \left( s^{\frac{1}{p_1}} \bar{f} (s)\right)^\sigma \frac{ds}{s}\bigg)^{\frac{1}{\sigma}} \right)^q \frac{d\gamma}{\gamma}  \right)^{\frac{1}{q}}.
$$

Further, for the first and third terms, we apply the following variants of Hardy inequalities: if  $\mu >0, -\infty < \nu < \infty$ and $0 <\sigma, \tau \le \infty$, then 
   
$$
   \left( \int_0^{\infty} \left( y^{-\mu} \left( \int_0^y \left(r^{-\nu} |g(r)| \right)^{\sigma} \frac{dr}{r} \right)^{\frac{1}{\sigma}} \right)^{\tau} \frac{dy}{y}\right) ^{\frac{1}{\tau}} \le \left(\mu \sigma \right)^{-\frac{1}{\sigma}} \left( \int_0^{\infty} \left( y^{-\mu-\nu} |g(y)| \right)^{\tau} \frac{dy}{y} \right) ^{\frac{1}{\tau}}
$$
   and
   $$
    \bigg( \int_0^{\infty} \bigg( y^{\mu} \bigg( \int_0^y \bigg(r^{-\nu} |g(r)| \bigg)^{\sigma} \frac{dr}{r} \bigg)^{\frac{1}{\sigma}} \bigg)^{\tau} \frac{dy}{y}\bigg) ^{\frac{1}{\tau}} \le \big(\mu \sigma \big)^{-\frac{1}{\sigma}} \bigg( \int_0^{\infty} \bigg( y^{\mu-\nu} |g(y)| \bigg)^{\tau} \frac{dy}{y} \bigg) ^{\frac{1}{\tau}}.
   $$

According to these inequalities we have
$$
\|f\|_{(N_{p_0, \sigma}(M), N_{p_1,\sigma}(M))_{\theta,q}} \lesssim \bigg( \int_0^{\infty} \bigg( \gamma^{\frac{1}{p}} \Bar{f}(\gamma) \bigg)^q \frac{d\gamma}{\gamma} \bigg)^{\frac{1}{q}} =  \|f\|_{N_{p,q}},
$$
and hence, 
$$
\|f\|_{(N_{p_0, q_0}(M),N_{p_1, q_1}(M))_{\theta,q}} \lesssim \|f\|_{(N_{p_0, \sigma}(M), N_{p_1, \sigma}(M))_{\theta,q}} \lesssim \|f\|_{N_{p,q}(M)},
$$
where $\frac{1}{p}= \frac{1-\theta}{p_0}+\frac{\theta}{p_1}.$ From the theorem in \cite[Theorem 1]{N1} we know that the following embedding holds
$$
(N_{p_0, q_0}(M), N_{p_1, q_1}(M))_{\theta,q} \hookrightarrow N_{p,q}(M),
$$
i.e. 
$$
\|f\|_{N_{p,q}(M)} \lesssim \|f\|_{(N_{p_0, q_0}(M),N_{p_1, q_1}(M))_{\theta,q}}.
$$

This proves the equivalence (\ref{th2}).

\end{proof}

\section{Proof of Corollary 1}
\begin{proof}

 According to the real interpolation method \cite[Theorem 3.1.2]{BL}  and the inequalities (\ref{th3eq1}) and (\ref{th3eq2}) it follows 
 $$
 \|Tf\|_{(N_{q_0,\infty}(M), N_{q_1,\infty}(M))_{\theta, \tau}} \le F^{1-\theta}_0 F^{\theta}_1 \|f\|_{(N_{p_0,\sigma}(M), N_{p_1,\sigma}(M))_{\theta, \tau}}.
 $$

 From the relation (\ref{eq1}) we have 
 $$
 \|Tf\|_{N_{q,\tau}(M)} \le c \|Tf\|_{(N_{q_0,\infty}(M), N_{q_1,\infty}(M))_{\theta, \tau}}.
 $$

 From Theorem \ref{th2}, taking into account, that $f\ge 0$, we obtain
 $$
 \|f\|_{N_{p,\tau}(M)} \asymp \|f\|_{(N_{p_0, \sigma}(M), N_{p_1, \sigma}(M))_{\theta,q}}.
 $$

\end{proof}

\end{document}